\documentclass[11pt]{article}

\usepackage{amssymb,amsmath,amsthm}

\usepackage{indentfirst}

\usepackage{amssymb,amsmath,amsthm}

\newtheorem{Th}{\hskip\parindent Theorem}
\newtheorem{Le}{\hskip\parindent Lemma}
\newtheorem{Zam}{\hskip\parindent Remark}
\newtheorem{Prop}{\hskip\parindent Proposition}

\newcommand{\M}{\mathcal M}
\newcommand{\Li}{{\rm Li}}
\renewcommand{\le}{\leqslant}\renewcommand{\ge}{\geqslant}
\author{V.Shur\footnote{Mathematics Department, Moscow State University,
Russia, vladimir@chg.ru }\and Ya.
Sinai\footnote{Mathematics Department, Princeton
University, USA, sinai@math.princeton.edu} \and A.
Ustinov\footnote{Khabarovsk Division of Institute for
Applied Mathematics, Far Eastern Branch of the Russian
Academy of Science, Russia, ustinov@iam.khv.ru }}

\title{Limiting Distribution of Frobenius Numbers for $n=3$}

\begin{document}

\maketitle \setcounter{Zam}0

\section{Introduction}

The purpose of this paper is to give a complete derivation
of the limiting distribution of large Frobenius numbers
outlined in~\cite{Bourgain2007} and fill some gaps
formulated there as hypotheses. We start with the basic
definitions and descriptions of some results.

\medskip
Consider $n$ mutually coprime positive integers $a_1$,
$a_2$, \dots, $a_n$. This means that there is no $r>1$ such
that each $a_j$, $1\le j\le n$, is divisible by $r$. Take
$N$ which later will tend to infinity and will be our main
large parameter. Introduce the ensemble $Q_N$ of mutually
coprime $a=(a_1,\dots,a_n)$, $1\le a_j\le N$, $1\le j\le n$
and $P_N$ be the uniform probability distribution on $Q_N$.
For each $a\in Q_N$ denote by $F(a)$ the largest integer
number that is not representable in the form
$x=x_1a_1+\cdots+x_na_n$, where $x_j$ are non-negative
integers. $F(a)$ can be considered as a random variable
defined on $Q_N$. The basic problem which will be discussed
in this paper is the existence and the form of the limiting
distribution for the normalized Frobenius numbers
$f(a)=\dfrac 1{N^{1+\frac{1}{n-1}}}F(a)$. The reason for this
normalization will be explained below.

\medskip
The case of $n=2$ is simple in view of the classical result
of Sylvester (see~\cite{Sylvester1884})
according to which
$F(a_1,a_2)=a_1a_2-a_1-a_2$. It shows that in a typical
situation $F$ grow as $N^2$. The first non-trivial case is
$n=3$ where $F(a)$ grow as $N^{3/2}$ It is known
(see~\cite{Ustinov2009}) that the numbers $F(a_1,a_2,a_3)$ have
weak asymptotics:
$$
  \frac 1{x_1x_2a_3^{7/2}}
  \sum_{a_1\le x_1a_3}
  \sum_{a_2\le x_2a_3}
  \left(F(a_1,a_2,a_3)-
  \dfrac 8\pi \sqrt{a_1a_2a_3}\right)=
  O_{x_1,x_2,\varepsilon}\left(a_3^{-1/6+\varepsilon}\right)
$$
For arbitrary $n$ the only result known to us is the
following theorem proven in~\cite{Bourgain2007}.

\medskip
\begin{Th}
\label{Th_1} Under some additional technical condition
(see~\cite{Bourgain2007}) the family of probability
distributions of $f_N(a)=\frac 1{N^{1+\frac 1{n-1}}}F(a)$
is weakly compact. This means that for every
$\varepsilon>0$ one can find $\mathcal D=\mathcal
D(\varepsilon)$ such that
$$
  P_N\left\{\dfrac1{N^{1+\frac 1{n-1}}}F(a)\le \mathcal D\right\}\ge 1-\varepsilon.
$$
\end{Th}

\medskip
In this theorem $\varepsilon,\mathcal D$ do not depend on
$N$. It also implies the existence of the limiting points
(in the sense of weak convergence) for the sequence of
probability distributions of $f_N(a)$. As was already
mentioned, in this paper we shall study the limiting
distribution of $f_N(a)=\frac 1{N^{3/2}}F(a)$,
$a=(a_1,a_2,a_3)$  as $N\to\infty$. This distribution is
not universal and will be described below.

\medskip
Take any $\rho$,  $0<\rho <1$, and consider its expansion
into continued fraction
\begin{equation}
\label{eq_1}
  \rho=[h_1,h_2,\dots,h_s,\dots]
\end{equation}
where $h_j\ge 1$ are integers. If $\rho$ is rational then
the continued fraction~\eqref{eq_1} is finite. The finite
continued fractions $\rho=[h_1,\dots,h_s]=\dfrac{p_s}{q_s}$
are called the $s$-approximants of $\rho$. The numbers
$q_s$ satisfy the recurrent relations
\begin{equation}
\label{eq_2} q_s=h_sq_{s-1}+q_{s-2},\;\;s\ge 2
\end{equation}
Introduce the Gauss measure on $[0,1]$ given by the density
$\pi(x)=\frac1{\ln 2(1+x)}$. Then the elements of the
continued fraction~\eqref{eq_1} become random variables. It
is well-known that their probability distributions are
stationary in the sense that the distributions of any
$h_{m-k}$,
$h_{m-k+1},\dots,h_m.\dots,h_{m+k}$ do not depend
on $m$. We shall need the values of $s =s_1$, such that
$q_{s_1}$ is the first $q_s$ greater than $\sqrt N$. It was
proven in~\cite{Sinai2008} that $q_{s_1}/\sqrt N$ have a
limiting distribution. More precisely, the following
theorem is true.

\medskip
\begin{Th}
\label{Th_2}
Let $k$ be fixed and $s(R)$ be the smallest $s$ for which
$q_s\ge R$. As $R\to\infty$ there exists the joint limiting probability
distribution of $\frac{q_{s(R)}}R$, $h_{s(R)-k}$, \dots,
$h_{s(R)+k}$.
\end{Th}

\medskip
In the paper~\cite{Ustinov2008b} the analytic form of this
distribution was given.

\medskip
Consider the subensemble $Q_N^{(0)}\subset Q_N$ for which
$a_1,a_3$ are coprime. Then there exists $a_1^{-1}
(\bmod\,a_3)$, $1\le a_1^{-1}<a_3$. Denote
$\rho=\frac{a_1^{-1}a_2}{a_3}$. The expansion of $\rho$
into continued fraction will be needed below. Clearly, $\rho$ is a
rational number. However, the following theorem is valid.

\medskip
\begin{Th}
\label{Th_3} As before, consider $s_1$ such that
$q_{s_1-1}<\sqrt N<q_{s_1}$. Then in the ensemble
$Q_N^{(0)}$ equipped with the uniform measure, for any
$k>0$ and $N\to\infty$ there exists the joint
limiting probability distributions of $\frac{q_{s_1}}{\sqrt
N}$, $h_{s_1-k}$, \dots, $h_{s_1+k}$ which coincides with
the distribution in theorem 2.
\end{Th}

\medskip
A stronger  version of theorem~\ref{Th_3} is also valid.

\medskip
\begin{Th}
\label{Th_4} Let the first elements of the continued
fraction for $\rho$ be fixed: $h_1,h_2,\dots,h_l$. Then
under this condition and as
$N\to\infty$ the conditional distributions of
$\frac{q_{s_1}}{\sqrt N}$, $h_{s_1-k}$, \dots, $h_{s_1+k}$
converge to the same limit as in theorems~\ref{Th_2}
and~\ref{Th_3}.
\end{Th}

\medskip
All these theorems will be proven in section~\ref{Sec_3}. Now
we can formulate the main result of this paper.

\medskip
\begin{Th}
\label{Th_5}There exists the limiting distribution of
$f_N(a)=f_N((a_1,a_2,a_3))$,
$(a_1 , a_2 , a_3) \in Q_N$ as $N
\rightarrow \infty$.\end{Th}

\medskip
The proof of the main theorem is given in
section~\ref{Sec_2}. First we consider the ensemble
$Q_N^{(0)}$ and then explain how to extend the proof to
$Q_N$.

\medskip
The second author thanks NSF for the financial support,
grant  DMS No 0600996. The research of the third author was
supported by the Russian Foundation for Basic Research (grant
no. 07-01-00306 
and the
Russian Science Support Foundation.

\medskip
\section{The limiting Distribution of $f_N(a)$.}\label{Sec_2}

\medskip
Return back to the case of arbitrary $n$. Introduce arithmetic
progressions
$$
  \Pi_r=\{r+ma_n,m\ge 0\},\quad 0\le r<a_n.
$$
For non-negative integers $x_1$, \dots, $x_{n-1}$ such that
$x_1a_1+x_2a_2+\cdots+x_{n-1}a_{n-1}\in\Pi_r$ we write
$$x_1a_1+\cdots+x_{n-1}a_{n-1}=r+m(x_1,\dots,x_{n-1})a_n.
$$
Define $\overline m(r)=\underset{x_1\dots,x_{n-1}}\min
m(x_1,\dots,x_{n-1})$ and put
$$
F_1(a)=\max_{0\le r<a_n} \quad
\min_{{x_1,\dots,x_{n-1}\atop
x_1a_1+\cdots+x_{n-1}a_{n-1}\in\prod_r}}
                                     (r+m(x_1,\dots,x_{n-1})a_n) =
$$
$$
= \, \max\limits_{0 \le r < a_n}
\quad \min\limits_{x_1 a_1 + \ldots + x_n + a_{n-1}
\:\equiv \: r
(\hspace{-.8em}\mod a_n)} (x_1 a_1 + \ldots + a_{n-1} a_{n-1} ) \, .
$$
It  was proven in~\cite{Selmer1978} that $F(a)=F_1(a)-a_n$.
A slightly weaker statement can be found
in~\cite{Bourgain2007}. Since in a typical situation $a_j$
grow as $N$ while  $F_1(a)$ grow as $N^{1+\frac 1{n-1}}$
(see also~\cite{Bourgain2007}) the limiting behavior of
$\frac {F(a)}{N^{1+\frac 1{n-1}}}$  and $\frac
{F_1(a)}{N^{1+\frac 1{n-1}}}$ is the same, but the analysis
of $\frac{F_1(a)}{N^{1+\frac 1{n-1}}}$ is slightly simpler.
Let us write for $n=3$
$$
  x_1a_1+x_2a_2=r+m(x_1,x_2)a_3
$$
or
\begin{equation}
\label{eq_3}
  x_1a_1+x_2a_2\equiv r(\bmod\,a_3)
\end{equation}
We assume that $a_1,a_3$ and $a_2, a_3$ are coprime. Therefore
there exists $a_1^{-1}$, $1\le a_1^{-1}<a_3$, such that
$a_1\cdot a_1^{-1}\equiv 1(\bmod\,a_3)$. Choose $a_1^{-1}$
so that $1\le a_1^{-1}<a_3$ and rewrite~\eqref{eq_3} as
follows
\begin{equation}
\label{eq_4}
  x_1+a_{12}x_2\equiv r_1(\bmod\,a_3)
\end{equation}
where $a_{12}\equiv a_1^{-1}a_2(\bmod\,a_3)$,
$0<a_{12}<a_3$ and $r_1\equiv ra_1^{-1}(\bmod\,a_3)$, $0\le
r_1<a_3$. From~\eqref{eq_4}
\begin{equation}
\label{eq_5}
  a_{12}x_2\equiv (r_1-x_1)(\bmod\,a_3)
\end{equation}
The expression~\eqref{eq_5} has a nice geometric
interpretation. Consider $S=[0,1,\dots,a_3-1]$ as a
``discrete circle''. Let $\mathcal R$ be the rotation of
this circle by $a_{12},$ i.e.\break $\mathcal Rx=x+a_{12}(\bmod
a_3)$. Then $\mathcal R^px=x+pa_{12}(\bmod a_3)$
and~\eqref{eq_5} means that $r_1-x_1$ belongs to the orbit
of $0$ under the action of $\mathcal R$. From the
definition of $F_1(a)$
\begin{gather}\nonumber
  F_1(a)=\max_{0\le r<a_3}\min_{{x_1a_1+x_2a_2\equiv
  r(\bmod\,a_3)\atop0\le x_1,x_2<a_3}}(x_1a_1+x_2a_2)=                                 \\\label{eq_6}
  =N^{3/2}\max_{0\le r_1<a_3}\min_{x_1+x_2a_{12}\equiv r_1\pmod{a_3}}
  \left(\dfrac{x_1}{\sqrt N}\dfrac{a_1}N+\dfrac{x_2}{\sqrt N}\dfrac{a_2}N\right)
\end{gather}
Choose $h^{(j)}=(h_1^{(j)},\dots,h_m^{(j)})$, $j=1,2,3$ and
denote by $Q_{N,h^{(1)},h^{(2)},h^{(3)}}^{(0)}$ the
ensemble of $a=(a_1,a_2,a_3)\in Q_N^{(0)}$ such that the
first $m$  elements of the continued fractions of
$\frac{a_j}N$  are given by $h^{(j)}$, $j=1,2,3$. This step
means the localization of the ensemble $Q_N^{(0)}$. It is
easy to see that for every $\varepsilon >0$  one can find
rational $\alpha_1$, $\alpha_2$, $\alpha_3$ and $m$ such
that $\left|\frac{a_j}N-\alpha_j\right|\le
\varepsilon,\;1\le j \le 3$. Then in~\eqref{eq_6} one can
replace $\frac{a_j}N$ by $\alpha_j$. Since
$\frac{x_j}{\sqrt N}$ will take the values $O(1)$ the whole
expression in~\eqref{eq_6} takes values $O(1)$ and instead
of~\eqref{eq_6} we may consider
\begin{equation}
\label{eq_7} \max_{r_1}\min_{x_1+a_{12}x_2\equiv
r_1\pmod{a_3}}
  \left(\dfrac{x_1}{\sqrt N}\alpha_1+\dfrac{x_2}{\sqrt N}\alpha_2\right)
\end{equation}
with the error $O(\varepsilon)$. We assume that in the ensemble
$Q_{N,h^{(1)},h^{(2)},h^{(3)}}^{(0)}$ we also have the
uniform distribution.

\medskip
We shall need some facts from the theory of rotations of
the circle. According to our assumption $a_{12}$ and $a_3$
are coprime. Therefore $\mathcal R$ is ergodic in the sense
that $\mathcal R^{a_3}=Id$ and $a_3$ is the smallest number
with this property. Put $\rho=\frac {a_{12}}{a_3}$  and
write down the expansion of $\rho$ into continued fraction:
$\rho=[h_1,h_2,\dots , h_{s_0}]$. Let also be
$\rho_s=[h_1,h_2,\dots,h_s]=\frac{p_s}{q_s}$ and $s_1$ is
such that $q_{s_1-1}<\sqrt N<q_{s_1}$.

\medskip
It will be more convenient to consider the usual unit
circle instead of $S$ and use the same letter $\mathcal R$
for the rotation of the unit circle by $\rho$. Introduce
the interval $\varDelta_0^{(p)}$ bounded by $0$ and
$\{q_p\rho\}$ and $\varDelta_j^{(p)}=\mathcal
R^j\varDelta_0^{(p)}$. Using the induction one can show
that $\varDelta_j^{(p)},\;0\le j<q_{p+1}$ and
$\varDelta_j^{(p+1)},\;0\le j'<q_p$ are pair-wise disjoint
and their union is the whole circle except the boundary
points (see~\cite{Sinai2008}). Denote by $\eta^{(p)}$ the
partition of the unit circle into
$\varDelta_j^{(p)},\;\varDelta_{j'}^{(p+1)}$. Then
$\eta^{(p+1)}\ge\eta^{(p)}$ in the sense that each element
of $\eta^{(p)}$ consists of several elements of
$\eta^{(p+1)}$. More precisely, $\varDelta_0^{(p-1)}$
consists of $h_p$ elements $\varDelta_j^{(p)}$ and one
elements $\varDelta_0^{(p+1)}$. The partitions $\eta^{(p)}$
show how the orbit of $0$ fills the circle.

\medskip
Return back to the discrete circle $S$. The partitions
$\eta^{(p)}$ can be constructed in the same way as in the continuous
case.
We have to analyze
\begin{equation}
\label{eq_8}
  \underset{0\le r_1<a_3}\max \quad
  \min_{\substack{x_1,x_2   \\
        x_1+a_{12}x_2\equiv r_1(\bmod\,a_3)}}
  \left(\dfrac{x_1}{\sqrt N}\alpha_1+\dfrac{x_2}{\sqrt N}\alpha_2\right)
\end{equation}
for given $\alpha_1,\alpha_2,\;0<\alpha_1,\alpha_2<1$.

\begin{Le}
\label{Le_1} There exists some number $C_1 (\alpha_1,
\alpha_2) = C_1$ such that for any $r_1$ the point $x_1$
giving $\min \left( \frac{x_1}{\sqrt{N}} \, \alpha_1 + \,
\frac{x_2}{\sqrt{N}} \, \alpha_2 \right)$ under the
condition is such that $r_1 - x_1$ ($x_1 + a_{12} x_2
\equiv r_1 \pmod{a_3}$) is an end-point of some element of
the partition $\eta^{(s_1 + m_1)}$.  Here $m_1 \geq 0$ is
such that ${q{s_1} + m_1}/{q_{s_1}} \, \leq \, C_1
(\alpha_1, \alpha_2 )$
\end{Le}

\medskip
The proof is simple.  In any case $r_1 - x_1$ is an end-point of some
element of the partition $\eta^{(s_1 + m_1)}$.  If $m_1$ is too big then
$\frac{x_2}{\sqrt{N}}$ is too big because it takes too much time to
reach an end-point of $\eta^{(s_1 + m_1)}$ which is not an end-point of
one of the previous partitions.  We can choose $y_1$ so that $r_1 - y_1$
will be an end-point of some element of $\eta^{(s_1)}$ and the linear
combination $\frac{y_1}{\sqrt{N}} \alpha_1 + \frac{y_2}{\sqrt{N}}
\alpha_2$ is smaller.  This completes the proof of the lemma.

\medskip
Its meaning is the following.  If $r_1 - x_1$ is an end-point of
$\eta^{(s_1 + m_1)}$ with too big $m_1$ then $x_2$ will be also too big.

\medskip
Lemma 2 shows that $x_1$  also cannot be too big.
\begin{Le}
\label{Le_2} There exists an integer $m_2 > 0$ depending on
$\alpha_1 , \alpha_2$ the ratio $q_{s_1}/N$
and the elements of the continued fraction $h_{s_1},
h_{s_1+1} , \ldots , h_{s_1 + m_2}$ of $\rho$ such that for
any $r_1$ the interval $[r_1 - x_1, r_1]$ corresponding to
the minimum of
\[
\frac{x_1}{\sqrt{N}} \, \alpha_1 \, + \, \frac{x_2}{\sqrt{N}} \,
\alpha_2
\]
has not more than $m_2$ elements of $\eta^{(s_1)}$.
\end{Le}

\medskip
The proof is also simple.  If the number in question is too big then
$\frac{x_1}{\sqrt{N}}$ will be too big.  Therefore for given $r_1$
$\min$ can be attained at a point which is closer to $r_1$.

\medskip
The values of $q_{s_1}/\sqrt{N}$
 and $h_{s_1} , h_{s_1 + 1} \ldots$,
$h_{s_1 + m_2}$ determine the structure of the
partitions\break $\eta^{(s_1)} , \ldots , \eta^{(s_1 +
m_2)}$.

\medskip
The conclusion which follows from both lemmas is that for each $r_1$ we
check only finitely many $x_1$ and $x_2$ and find $\min ( x_1 \alpha_1 +
x_2 \alpha_2)$ among them.  The number of points which have to be
checked depends on $\alpha_1, \alpha_2$, $\frac{q_{s_1}}{\sqrt{N}}$ and
$h_{s_1} , \ldots , h_{s_1 + m_2}$.

\medskip
Now we remark that $r_1$ must be also an end-point of some element of
the partition $\eta^{(s_1)}$.  Indeed, if $r_1$ increases within some
element of $\eta^{(s_1)}$ then the set of values $r_1 - x_1$ which have
to be checked remain the same.  Then $\max\limits_{r_1}$ is attained at
the end-point of this element $\eta^{(s_1)}$ because  $r_1 - x_1$ is a
monotone increasing function of $r_1$.

\medskip
The last step in the proof is the final choice of $r_1$.
As was mentioned above $r_1$ must be an end-point of some
element of $\eta^{(s_1)}$ and $\frac{x_1}{\sqrt{N}}$ takes
finitely many values. Therefore $r_1$ should be chosen so
that $x_2/{\sqrt{N}}$ takes the largest possible value.
Take the last point $r^\prime_1 = \mathcal{R}^{q_{s_1 - 1}}
0$ on the orbit of $0$ of the length $q_{s_1}$. Assume for
definiteness that $r^\prime_1$ lies to the left from $0$.
Consider $m_2$ elements of $\eta^{(s_1)}$ which start from
$r^\prime_1$ and go left.  Then $r_1$ must be one of the
end-points of these elements.  Indeed, if $r_1$ lies more
to the left from $0$ then the values $x_1$ take finitely
many values and $x_2$ will be significantly smaller.
Therefore it cannot give maximum over $r$ of our basic
linear form.

\medskip
Thus we take $m_2$ elements of $\eta^{(s_1)}$, consider their
end-points.
Each end-point is a possible value of $r$.  Taking finitely many $x_1$
(see Lemma 1 and Lemma 2) we find minimum of our basic linear form.
After that we find $r$ for which this minimum takes maximal value.  In
this way we get the solution of our max-min problem.  It is clear that
this solution is a function of $\frac{q_{s_1}}{\sqrt{N}}$ and elements
$h_j , s_1 \leq j \leq s_1 + m_1$ of the continued fraction of $\rho$
near $s_1$.  Since $\frac{q_{s_1}}{\sqrt{N}}$ and $h_j , s_1 \leq j \leq
s_1 + m_1$ have limiting distribution as $N \rightarrow \infty$ the
number $f_N ( a ) = \frac{1}{N^{3/2}} \, F_1 ( a )$ also has a limiting
distribution.

\medskip
It remains to extend our proof to the case when the pairs from
$a_1 , a_2 , a_3$ have
non-trivial common divisors, say $k_1$ is $gcd$ of $a_1 , a_3$ and $k_2$
is $gcd$ of $a_2, a_3$. It is easy to show that $k_1 , k_2$ have a joint
limiting probability distribution in the whole ensemble $Q_N$.
Fixing $k_1 , k_2$ we can write $a_1
= k_1 a^\prime_1$, $a_2 = k_2 a^\prime_2$, $a_3 = k_1 k_2 a^\prime_3$
where $a^\prime_1 , a^\prime_3$ are coprime, $a^\prime_2 , a_3$ are
coprime and $k_1 , k_2$ are coprime. This implies that
$(a^\prime_1)^{-1}$$(\hspace{-.3em}\mod a^\prime_3)$ exists and we can
multiply both sides of (3) by $(a^\prime_1)^{-1}$.  This will give
\begin{equation}
k_1 x_1 + k_2 a^\prime_2 \cdot (a^\prime_1)^{-1} \cdot x_2 \, \equiv \,
r_1 (\hspace{-.8em}\mod a_3)
\end{equation}
where $r_1 = r \cdot ( a^\prime_1)^{-1}$ $( \hspace{-.5em}\mod a_3)$.
Denote $b = a^\prime_2 ( a^\prime_1)^{-1}$.

\medskip
Then from (9) we have the linear form
\begin{equation}
k_1 x_1 + k_2 \, b x_2 \, \equiv \, r_1 \, ( \hspace{-.8em}\mod a_3)
\end{equation}
which we can treat in the same way as before.
\section{Statistical properties of continued
fractions}\label{Sec_3}

\medskip
Statistical properties of elements of continued fractions
usually are identical for real numbers and for rationales
with bounded denominators
(see~\cite{Ustinov2005a}--\cite{Ustinov2008a}).

\medskip
Let $\M$ be the set of integer matrices $S=\bigl(\begin{smallmatrix}
P & P'\\
Q & Q'
\end{smallmatrix}\bigr)$
with determinant $\det S=\pm 1$ such that
$ 1\le Q\le Q'$, $ 0\le P\le Q$, $ 1\le P'\le Q'.$
For real $\alpha\in(0,1)$ the fractions $P/Q$ and $P'/Q'$ with
$S=\bigl(\begin{smallmatrix}
P & P'\\
Q & Q'
\end{smallmatrix}\bigr)\in\M$ will be consecutive
convergents  to $\alpha$ (distinct from $\alpha$) if and only if
$$
0<\frac{Q'\alpha-P'}{-Q\alpha+P}= S^{-1}(\alpha)<1
$$ (see \cite[lemma 1]{Ustinov2005a}). Moreover if $\alpha=[0;h_1,h_2,\ldots]$ then for some $s\ge 1$
\begin{align}\label{MA}
\frac PQ=[0;h_1,\dots,h_{s-1}],&\quad \frac
{P'}{Q'}=[0;h_1,\dots,h_s],\\\dfrac{Q}{Q'}=[0;h_s, \ldots,
h_1],&\nonumber\quad
\frac{Q'\alpha-P'}{-Q\alpha+P}=[0;h_{s+1}, h_{s+2},
\ldots].
\end{align}
It means that the distribution of partial quotients  $h_{s-k}$,
\ldots, $h_{s+k}$ depends on Gauss-Kuz'min statistics of
fractions $Q/Q'$ and $({Q'\alpha-P'})/({-Q\alpha+P})$.

\medskip
For real $\alpha$, $x_1$, $x_2$, $y_1$, $y_2\in(0,1)$
denote by $N_{x_1,x_2,y_1,y_2}(\alpha,R)$ the number of
solutions of the following system of inequalities
\begin{gather}
0<S^{-1}(\alpha)\le x_1,\label{31_1}\quad Q\le x_2Q',\quad Q\le y_1R,\quad R\le
y_2Q',
\end{gather} with variables $P$, $P'$, $Q$, $Q'$ such that $S=\bigl(\begin{smallmatrix}
P & P'\\
Q & Q'
\end{smallmatrix}\bigr)\in\M$. Let
$$N(R)=N_{x_1,x_2,y_1,y_2}(R)=\int_0^1 N_{x_1,x_2,y_1,y_2}(\alpha,R)\,d\alpha$$
and
$$G(x_1,x_2,y_1,y_2)=\begin{cases}\frac{2}{\zeta(2)}\left(\log(1+x_1x_2)\log\frac{y_1y_2}{x_2}-\Li_2(-x_1x_2)\right),&\text{if }x_2\le y_1y_2;\\
-\frac{2}{\zeta(2)}\Li_2(-x_1y_1y_2),&\text{if }x_2> y_1y_2,
\end{cases}$$
where $\Li_2(\cdot)$ is the dilogarithm
$$\Li_2(z)=\sum\limits_{k=1}^{\infty}\dfrac{z^k}{k^2}=
-\int_{0}^{z}\dfrac{\log(1-t)}{t}dt.$$

\medskip
The next statement implies Theorem~\ref{Th_2}.

\begin{Prop}
\label{Le031.1} For $R\ge2$
$$N(R)=G(x_1,x_2,y_1,y_2)+O\left(\dfrac{x_1\log R}{R}\right).$$\end{Prop}

\medskip
\begin{proof} For every number $\alpha=[0;a_1,a_2,\ldots]$
find a unique matrix $S\in\M$ with elements $P$, $P'$, $Q$, $Q'$ defined
by~\eqref{MA} with the additional restriction $Q\le R<Q'$. The inequalities
$0<S^{-1}(\alpha)\le x_1$ define the interval $I_{x_1}(S)\subset(0,1)$ of the
length
$$|I_{x_1}(S)|=\left|\frac{P'+x_1P}{Q'+x_1Q}-\frac{P'}{Q'}\right|=
\frac{x_1}{Q'(Q'+x_1Q)}.$$ Hence
$$N(R)=\sum\limits_{\bigl(\begin{smallmatrix}
P & P'\\
Q & Q'
\end{smallmatrix}\bigr)\in\M}[Q\le x_2Q',Q\le y_1R,R\le y_2Q']
\frac{x_1}{Q'(Q'+x_1Q)},$$ where $[A]$ is $1$ if the
statement $A$ is true, and it is $0$ otherwise. Second row
$(Q,Q')$ can be complemented to the matrix from $\M$ in two
ways. That is why
\begin{equation}\label{Art_31.3}
N(R)=2\sum\limits_{Q'\ge R/y_2}\sum\limits_{{(Q,Q')=1}}[Q\le x_2Q',Q\le
y_1R]\frac{x_1}{Q'(Q'+x_1Q)}.
\end{equation}
In the first case $x_2\le y_1y_2$ and the M\"obius inversion formula gives
\begin{align*}N(R)=&2\sum\limits_{d\le
R}\dfrac{\mu(d)}{d^2}\sum\limits_{R/(y_2d)\le Q'< y_1R/(x_2d)}\sum\limits_{Q\le
x_2Q'}\frac{x_1}{Q'(Q'+x_1Q)}+\\+&2\sum\limits_{d\le
R}\dfrac{\mu(d)}{d^2}\sum\limits_{Q'\ge y_1R/(x_2d) }\sum\limits_{Q\le
y_1R/d}\frac{x_1}{Q'(Q'+x_1Q)}=\\=&
\frac{2}{\zeta(2)}\left(\log(1+x_1x_2)\log\frac{y_1y_2}{x_2}+
\int_{1/(x_1x_2)}^{\infty}\log\left(1+\frac1t\right)\frac{dt}t\right)+
O\left(\dfrac{x_1\log R}{R}\right)=\\=&
\frac{2}{\zeta(2)}\left(\log(1+x_1x_2)\log\frac{y_1y_2}{x_2}-\Li_2(-x_1x_2)\right)+
O\left(\dfrac{x_1\log R}{R}\right).
\end{align*}
The second case $x_2> y_1y_2$ can be treated in the same way.
\end{proof}

\medskip
Let
$$L(R)=L_{x_1,x_2,y_1,y_2}(R)=\sum\limits_{b\le
R^2}\sum\limits_{{a\le
b\atop(a,b)=1}}N_{x_1,x_2,y_1,y_2}\left(\frac
ab,R\right).$$ Theorem~\ref{Th_3} will be proved in the
following form.

\begin{Prop}
\label{Le031.2} For $R\ge2$
$$\frac{2\zeta(2)}{R^4}L(R)=G(x_1,x_2,y_1,y_2)+O\left(\dfrac{x_1\log^2
R}{R}\right).$$\end{Prop}

\begin{proof} Let $\alpha=a/b$ be a given number and $S=\bigl(\begin{smallmatrix}
P & P'\\
Q & Q'
\end{smallmatrix}\bigr)\in\M$ be a solution of the system~\eqref{31_1}.
Denote by $m$ and $n$
the integers such that $mP+nP'=a, mQ+nQ'=b$. Then the system~\eqref{31_1} can be
written as follows
\begin{gather*}
mP+nP'=a,\quad mQ+nQ'=b,\\ 0<m/n\le x_1,\quad 0<Q/Q'\le x_2,\quad Q\le
y_1R,\quad R\le y_2Q'.
\end{gather*}
Summing up solutions of this system over $a$ and $b$ we get that the sum $L(R)$
equals to the number of solutions of the following system
\begin{gather*}
mQ+nQ'\le R^2,\quad 0<m/n\le x_1,\quad 0<Q/Q'\le x_2,\quad Q/y_1\le R<y_2Q',
\end{gather*}
where $\bigl(\begin{smallmatrix}
P & P'\\
Q & Q'
\end{smallmatrix}\bigr)\in\M$, $0\le m\le n$, $(m,n)=1$.
For given $Q$ and $Q'$ values of $P$ and $P'$ can be founded in two ways.
Number of solutions of the last system is equal to the area of the corresponding
region with the factor $1/\zeta(2)$ (see~\cite[Ch.~II,
problems~21--22]{Vinogradov1972})
$$\frac{R^4}{2\zeta(2)}\cdot\frac{x_1}{Q'(Q'+x_1Q)}+O\left(\dfrac{x_1
R^{2}\log R}{Q'}\right).$$ It leads to the sum similar to ~\eqref{Art_31.3}:
\begin{align*}L(R)=
\frac{R^4}{\zeta(2)}\sum\limits_{R/y_2\le Q'\le
R^2}\sum\limits_{{Q\le\min\{y_1R,x_2Q'\}\atop(Q,Q')=1}}\frac{x_1}{Q'(Q'+x_1Q)}+
O(x_1R^{3}\log^2R).
\end{align*}
Therefore
\begin{align*}L(R)=
\frac{R^4}{\zeta(2)}N(R)+ O(x_1R^{3}\log^2R),
\end{align*}
and Proposition~\ref{Le031.2} follows from Proposition~\ref{Le031.1}.
\end{proof}

In order to prove theorem 4 we have to use Kloosterman sums
$$K_q(m,n)=\sum\limits_{x,y=1}^{q}\delta_q(xy-1)\,e^{2\pi
i\frac{mx+ny}{q}},$$ where $\delta_q(a)$ is characteristic
function of divisibility by $q$:
$$\delta_q(a)=[a\equiv 0\!\!\pmod{p}]=\begin{cases}
1,& \text{if } a\equiv 0\pmod{q},\\
0,& \text{if } a\not\equiv 0\pmod{q}.\\
\end{cases}$$
Using Estermann bound (see~\cite{Estermann1961})
\begin{equation*}
\label{Estermannn} |K_q(m,n)|\le\sigma_0(q)\cdot(m,n,q)^{1/2}\cdot q^{1/2}.
\end{equation*} it is easy to prove the following statement
(see~\cite{Ustinov2008a} for details).

\begin{Le}
\label{LeKloo} Let $q\ge1$ be an integer, $Q_1$, $Q_2$, $P_1$, $P_2$ be real
numbers and $0\le P_1,P_2\le q$. Then the sum
$$\Phi_q(Q_1,Q_2;P_1,P_2)=\sum\limits_{
Q_1<u\le Q_1+P_1\atop Q_2<v\le Q_2+P_2 }\delta_q(uv-1)$$ satisfies the
asymptotic formula
$$\Phi_q(Q_1,Q_2;P_1,P_2)=\dfrac{\varphi(q)}{q^2}\cdot P_1P_2+O\left(\psi(q)\right),$$
where
\begin{equation*}
\label{Art_24.psi} \psi(q)=\sigma_0 (q)\log^2 (q+1)q^{1/2}.
\end{equation*}
\end{Le}

\medskip
It implies the following general result (see~\cite{Ustinov2005a}).

\begin{Le}
\label{Le3.2} Let $q\ge1$ be an integer and let  $a(u,v)$ be a function
defined on the set of integral points $(u,v)$ such that $1\le u,v\le q$. Assume that this
function satisfies the inequalities
\begin{equation}
\label{3.1}a(u,v)\ge0,\quad \Delta_{1,0}a(u,v)\le0,\quad
\Delta_{0,1}a(u,v)\le0,\quad \Delta_{1,1}a(u,v)\ge0
\end{equation} at all points at which these conditions have the
well-defined meaning. Then the
sum $$W=\sum\limits_{u,v=1}^{q}\delta_q(uv-1)a(u,v)$$ satisfies the
asymptotics
$$W=\dfrac{\varphi(q)}{q^2}\sum\limits_{u,v=1}^{q}a(u,v)+O\left(A\psi(q)\sqrt{q}\right),$$
where $\psi(q)$ is the function from lemma~\ref{LeKloo} and $A=a(1,1)$ is the
maximum of the function $a(u,v)$.
\end{Le}

Let
\begin{align*}
N_z(R)=&N_{z,x_1,x_2,y_1,y_2}(R)=\int_0^z
N_{x_1,x_2,y_1,y_2}(\alpha,R)\,d\alpha,\\
L_z(R)=&L_{z,x_1,x_2,y_1,y_2}(R)=\sum\limits_{b\le R^2}\sum\limits_{{a\le
zb\atop(a,b)=1}}N_{x_1,x_2,y_1,y_2}\left(\frac ab,R\right).
\end{align*}

The next statement implies Theorem~\ref{Th_4}.

\begin{Prop}
\label{Le031.3} For $R\ge2$ \begin{align*} N_z(R)=&z\cdot
G(x_1,x_2,y_1,y_2)+O\left(\dfrac{x_1\log^3R}{R^{1/2}}\right),\\
\frac{2\zeta(2)}{R^4}L_z(R)=&z\cdot
G(x_1,x_2,y_1,y_2)+O\left(\dfrac{x_1\log^3R}{R^{1/2}}\right).
\end{align*}\end{Prop}

\begin{proof}
Let $$\M_z=\left\{\begin{pmatrix}
P & P'\\
Q & Q'
\end{pmatrix}\in\M:\dfrac{P'}{Q'}\le z\right\}.$$ For a given $z$
there is at most one matrix $S=\bigl(\begin{smallmatrix}
P & P'\\
Q & Q'
\end{smallmatrix}\bigr)\in\M$ such that $Q\le R<Q'$ and $z\in I_{x_1}(S)$. Hence
\begin{align*}
N_z(R)=&\sum\limits_{\bigl(\begin{smallmatrix}
P & P'\\
Q & Q'
\end{smallmatrix}\bigr)\in\M_z}[Q\le x_2Q',Q\le y_1R,R\le y_2Q']
\frac{x_1}{Q'(Q'+x_1Q)}+O\left(\dfrac{x_1}{R^2}\right).
\end{align*}
If $Q'$ is fixed then $P'$ and $Q$ satisfy the congruence $P'Q\equiv
\pm1\pmod{Q'}$. Therefore
\begin{align*}
N_z(R)=&\sum\limits_{Q'\ge
R/y_2}\sum\limits_{P',Q=1}^{Q'}\delta_{Q'}(P'Q\pm1)[Q\le
\min\{x_2Q',y_1R\},P'\le zQ']
\frac{x_1}{Q'(Q'+x_1Q)}+O\left(\dfrac{x_1}{R^2}\right).
\end{align*}
Using Lemma~\ref{Le3.2} we obtain
\begin{align*}
N_z(R)=&\sum\limits_{Q'\ge
R/y_2}\dfrac{\varphi(Q')}{(Q')^2}\sum\limits_{P',Q=1}^{Q'}[Q\le
\min\{x_2Q',y_1R\},P'\le zQ']
\frac{x_1}{Q'(Q'+x_1Q)}+O\left(\dfrac{x_1\log^3R}{R^{1/2}}\right)=\\=&
z\sum\limits_{Q'\ge R/y_2}\dfrac{\varphi(Q')}{Q'}\sum\limits_{Q=1}^{Q'}[Q\le
\min\{x_2Q',y_1R\}]
\frac{x_1}{Q'(Q'+x_1Q)}+O\left(\dfrac{x_1\log^3R}{R^{1/2}}\right).
\end{align*}
Applying the formula
\begin{equation} \label{Art_31.Euler}
\dfrac{\varphi(Q')}{Q'}=\sum\limits_{d\mid Q'}\dfrac{\mu(d)}{d}
\end{equation} we get the same sum as in the proof of
Proposition~\ref{Le031.1}.

As in Proposition~\ref{Le031.2} the sum $L_z(R)$ equals to the number of
solutions of the system
\begin{gather*}
mQ+nQ'\le R^2,\quad mP+nP'\le z(mQ+nQ'),\\
0<m/n\le x_1,\quad 0<Q/Q'\le x_2,\quad Q/y_1\le R<y_2Q',
\end{gather*}
where $\bigl(\begin{smallmatrix}
P & P'\\
Q & Q'
\end{smallmatrix}\bigr)\in\M$, $0\le m\le n$, $(m,n)=1$. Again, there is at
most one matrix $S=\bigl(\begin{smallmatrix}
P & P'\\
Q & Q'
\end{smallmatrix}\bigr)\in\M$ such that $Q\le R<Q'$ and $z\in I_{x_1}(S)$. Also
for $Q'\ge R$
$$\sum\limits_{n\ge 1}\sum\limits_{m\le x_1n}[mQ+nQ'\le R^2]\ll x_1R^2.$$ This
estimate implies that
\begin{align*}L_z(R)=&
\frac{R^4}{\zeta(2)}\sum\limits_{\bigl(\begin{smallmatrix}
P & P'\\
Q & Q'
\end{smallmatrix}\bigr)\in\M_z}[R/y_2\le Q'\le R^2,Q\le\min\{y_1R,x_2Q'\}]
\frac{x_1}{Q'(Q'+x_1Q)}+ O(x_1R^{3}\log^2R)=\\=& \frac{R^4}{\zeta(2)}
\sum\limits_{R/y_2\le Q'\le
R^2}\sum\limits_{P',Q=1}^{Q'}[Q\le\min\{y_1R,x_2Q'\},P'\le
zQ']\frac{x_1\delta_{Q'}(P'Q\pm1)}{Q'(Q'+x_1Q)}+ O(x_1R^{3}\log^2R).
\end{align*}
Using Lemma~\ref{Le3.2} one more time we obtain

\begin{align*}
L_z(R)  = & \frac{R^4}{\zeta(2)}\sum\limits_{Q'\ge
R/y_2}\dfrac{\varphi(Q')}{(Q')^2}\sum\limits_{P',Q=1}^{Q'}[Q\le
\min\{x_2Q',y_1R\},P'\le zQ'] \frac{x_1}{Q'(Q'+x_1Q)}  + O
\left(x_1R^{7/2}\log^3R\right) =\\ =&
\dfrac{zR^4}{\zeta(2)}\sum\limits_{Q'\ge
R/y_2}\dfrac{\varphi(Q')}{Q'}\sum\limits_{Q=1}^{Q'}[Q\le
\min\{x_2Q',y_1R\}] \dfrac{x_1}{Q'(Q'+x_1Q)}  +
O\left(x_1R^{7/2}\log^3R\right). \end{align*}

\medskip
\noindent
Applying formula~\eqref{Art_31.Euler} we get the same sum as in as in
the proof of Proposition~\ref{Le031.1}.
\end{proof}

\begin{Zam}
In the simplest case  $x_2=y_1=y_2=1$ we have cumulative distribution function
$$F(x)=F(x,1,1,1)=-\frac{2}{\zeta(2)}\Li_2(-x),$$
which is not equal to the Gaussian function $\log_2(1+x)$.
As $x\to 0$ function $F(x)$ (with error terms in
Propositions~\ref{Le031.1} and ~\ref{Le031.2}) decreases as
a linear function $F(x)\sim 2x/\zeta(2)$. This fact shows
that the expectation of the partial quotient $a_s$ (defined
by inequalities $q_{s-1}\le R< q_s$) is equal to
infinity.\end{Zam}

\section{Concluding remarks}

The calculations done by one of the authors (A. Ustinov) shows that the
density of the limiting distribution of $\frac{F(a_1, a_2 ,
a_3}{\sqrt{a_1 a_2 a_3}}$ has the following simple form:
\[
p ( t ) =  \left\{
\begin{array}{ll}
0 , & \mbox{if} \ t \in [ 0, \sqrt{3}]; \\
\frac{12}{\pi} \, \left( \frac{t}{\sqrt{3}} - \sqrt{4-t^2} \right) , &
\mbox{if} \ t \in [ \sqrt{3} , 2]; \\
\frac{12}{\pi^2} \left( t \sqrt{3} \mbox{arccos} \, \frac{t
+ 3 \sqrt{t^2 - 4}}{4 \sqrt{t^2 - 3}}  +  \frac{3}{2}
\sqrt{t^2 - 4} \log \, \frac{t^2-4}{t^2 - 3} \right) , &
\mbox{if} \ t \in [ 2 , + \infty ) .
\end{array}
\right.
\]
This result will be published elsewhere.

\bigskip

\vfill
\today \ :gpp
\end{document}